\documentclass[11pt]{amsart}
\usepackage{amssymb,amsmath, hyperref, amsthm} 
\usepackage[normalem]{ulem} 

\usepackage{tikz}
\usepackage{ifthen}
\usepackage{enumitem} 

\newtheorem{theorem}{Theorem}

\newtheorem{conjecture}[theorem]{Conjecture}
\newtheorem*{conjecture*}{Conjecture}
\newtheorem{corollary}[theorem]{Corollary}
\newtheorem{definition}[theorem]{Definition}
\newtheorem{example}[theorem]{Example}

\newtheorem{exercise}[theorem]{Exercise}
\newtheorem{fact}[theorem]{Fact}
\newtheorem{lemma}[theorem]{Lemma}
\newtheorem{problem}[theorem]{Problem}
\newtheorem{proposition}[theorem]{Proposition}
\newtheorem{question}[theorem]{Question}
\newtheorem{remark}[theorem]{Remark}

\newcommand{\bcon}{\begin{conjecture}}
\newcommand{\econ}{\end{conjecture}}
\newcommand{\bcor}{\begin{corollary}}
\newcommand{\ecor}{\end{corollary}}
\newcommand{\bdf}{\begin{definition}}
\newcommand{\edf}{\end{definition}}
\newcommand{\benu}{\begin{enumerate}}
\newcommand{\eenu}{\end{enumerate}}
\newcommand{\beq}{\begin{equation}}
\newcommand{\eeq}{\end{equation}}
\newcommand{\bexa}{\begin{example}}
\newcommand{\eexa}{\end{example}}
\newcommand{\bexe}{\begin{exercise}}
\newcommand{\eexe}{\end{exercise}}
\newcommand{\bfac}{\begin{fact}}
\newcommand{\efac}{\end{fact}}
\newcommand{\bite}{\begin{itemize}}
\newcommand{\eite}{\end{itemize}}
\newcommand{\blem}{\begin{lemma}}
\newcommand{\elem}{\end{lemma}}
\newcommand{\bmat}{\begin{matrix}}
\newcommand{\emat}{\end{matrix}}
\newcommand{\bprb}{\begin{problem}}
\newcommand{\eprb}{\end{problem}}
\newcommand{\bpro}{\begin{proposition}}
\newcommand{\epro}{\end{proposition}}

\newcommand{\bque}{\begin{question}}
\newcommand{\eque}{\end{question}}
\newcommand{\brem}{\begin{remark}}
\newcommand{\erem}{\end{remark}}
\newcommand{\bthm}{\begin{theorem}}
\newcommand{\ethm}{\end{theorem}}

\newcommand{\bpr}{\begin{proof}}
\newcommand{\epr}{\end{proof}}

\newcommand{\comment}[1]{\,}

\newcommand{\Z}{\mathbb Z}

\newcommand{\R}{\mathbb R}

\newcommand{\ve}{\varepsilon}

\title{Fermat-Catalan and Tijdeman-Zagier conjectures for products}

\author{Adam S. Sikora}

\begin{document}

\thispagestyle{empty}

\begin{abstract}
We propose conjectural generalizations of the Fermat-Catalan conjecture, the Tijdeman-Zagier conjecture, and of the Fermat Last Theorem, in which powers are replaced by products of integers. We also formulate a new explicit version of the abc conjecture.
\end{abstract}

\address{244 Math Bldg, University at Buffalo, SUNY, Buffalo, NY 14260}
\email{asikora@buffalo.edu}

\pagestyle{myheadings}

\maketitle

%
\section{Introduction}
%

Fundamental to number theory is the understanding of the interplay between multiplications and additions of integers.
This interplay not only underpins basic arithmetic but also leads to deep mathematical insights and challenges, especially when considering sums of powers of integers.
Exploration of sums of powers has historically led to some of the most formidable problems in mathematics, 
such as the Fermat's Last Theorem, \cite{Wi}, the Fermat-Catalan, Tijdeman-Zagier, and Pillai's conjectures,  cf. \cite{DG, El,Na}. In this paper we propose conjectural generalizations of these statements to products, $X=x_1\cdot ... \cdot x_n,$
with $max(x_1,...,x_n)-min(x_1,...,x_n)$, called spread, bounded from above. 

Additionally, we propose a new explicit version of the abc conjecture (Conjecture \ref{c.abc-mine})
and discuss the implications of the abc conjectures to our work.
Specifically, we show that our Fermat-Catalan Conjecture for Products is a consequence the abc conjecture, see Theorem \ref{t.FCabc}. However, we do not expect that to be the case for our generalization of the Tijdeman-Zagier conjecture (Conjectures \ref{c.GBTZ}), nor for our Conjectures \ref{c.nonmaxgcd3}, \ref{c.Fp}, and \ref{c.main}.
Before we formulate these statements, we recall the relevant classical conjectures. 
\vspace*{.1in}

{\bf Acknowledgements.} The author would like to thank D. Zagier for a helpful discussion. 
The author was supported by the Simons Foundation grant 957582.
Computational testing of the conjectures of this paper was partially performed on the cluster computers at the Center for Computational Research at University at Buffalo.

%
\section{Fermat-Catalan and Tijdeman-Zagier conjectures}
%

The Fermat-Catalan conjecture asserts that 
\beq\label{e.FC}
x^{n}+y^{m}=z^{k}
\eeq
has only finitely many solutions  $(x^{n},y^{m},z^{k})$,  
with coprime $x,y,z$ and with $n,m,k$ such that 
\beq\label{e.suminverse}
\frac1{n}+\frac1{m}+\frac1{k}<1.
\eeq
In this paper, all letters, such as, $x,y,z,n,m,k$, (possibly with indices) indicate positive integers.

The known coprime solutions of the Fermat-Catalan equation \eqref{e.FC} (subject to \eqref{e.suminverse}) are: 
$$1+2^{3}=3^{2}, \quad 2^{5}+7^{2}=3^{4}, \quad  7^{3}+13^{2}=2^{9}, \quad 2^{7}+17^{3}=71^ {2}, \quad  3^{5}+11^{4}=122^{2},$$  
$$33^{8}+1549034^{2}=15613^{3}, \quad  1414^{3}+2213459^{2}=65^{7}, \quad  9262^{3}+15312283^{2}=113^{7},$$ 
$$17^{7}+76271^{3}=21063928^{2}, \quad 43^{8}+96222^{3}=30042907^{2}.$$
(Darmon-Granville,  \cite{DG}, proved that  \eqref{e.FC} has finitely many coprime solutions for any fixed $n,m,$ and $k$ satisfying \eqref{e.suminverse}.)

Based on the testing performed by us (see below) and others before us, it is reasonable to believe that the above are the only coprime solutions of \eqref{e.FC} subject to \eqref{e.suminverse}. 

The Tijdeman-Zagier conjecture, also known as Beal conjecture, states that \eqref{e.FC} has no coprime solutions for $n,m,k>2,$
\cite{El}. 
Andrew Beal offered a \$1,000,000 cash award for a resolution of this conjecture.

A more general version of Fermat-Catalan conjecture posits that for any positive $A,B,$ and $C$,
the equation 
\beq\label{e.FCgen}
A\cdot x^{n}+B\cdot y^{m}=C\cdot z^{k}
\eeq
has only finitely many coprime solutions  $(x^{n},y^{m},z^{k})$,  
with $n,m,k$ satisfying \eqref{e.suminverse}.
\vspace*{.1in}

\noindent{\bf Computer Verification.} 
We have verified that the above listed are the only coprime solutions of \eqref{e.FC} subject to \eqref{e.suminverse}
(and, hence, in particular, we verified the Tijdeman-Zagier  conjecture) when $min(n,m,k)\leq 113$ for all
$$x^{n},y^{m},z^{k}<M_{min(n,m,k)},$$
where 
$$M_2=2^{71},\  M_3=2^{80},\ M_4=2^{100}, \ \text{and}\ M_d=2^{113}, \ \text{for}\ 5\leq d\leq 113.$$ 

For comparison, the largest power appearing in the above solutions, $30042907^{2}$, is approximately $2^{49.7}.$
(We did not find any information in the existing literature on the scope of the earlier verifications of the Tijdeman-Zagier conjecture nor the earlier searches of coprime solutions of  \eqref{e.FC} subject to \eqref{e.suminverse}.)
\vspace*{.1in}

\section{The Fermat-Catalan and Pillai's conjectures for products}
\label{s.FCp}

An expression $X=x_1\cdot .... \cdot x_n$ (with $x_1,...,x_n\in \Z_{> 0}$, by our convention)
is a product of \underline{degree} $d_X=n$. The difference, 
$$s_X=max(x_1,...,x_n)-min(x_1,...,x_n),$$ 
is called its \underline{spread}. We will call $min(x_1,...,x_n)$ the \underline{base} of $X$ and denote it by $b_X.$

In particular, powers are products of spread zero. 
On the other hand, $6$ can be written as a product $2\cdot 3$ of degree $2$ and of spread $1$ or as a  product $1^k\cdot 2\cdot 3$ of any degree $\geq 3$ and of spread $2$.  

We conjecture that products of relatively small spread have certain properties of integral powers. 
In particular, we propose:

\bcon[Fermat-Catalan Conjecture for Products]\ \\ \label{c.FCS}
For any positive $A,B,C\in \Z$ and $F, Q,M\in \R$, $F<1,$  the equation 
$$A\cdot X+B\cdot Y=C\cdot Z$$ 
has finitely many solutions only,
for products $X,Y,Z$ of degrees $n_X,n_Y,n_Z$ respectively,  satisfying 
\beq\label{e.sBound}
\frac{s_X^2}{b_X},\ \frac{s_Y^2}{b_Y},\ \frac{s_Z^2}{b_Z} \leq M,
\eeq 
\beq\label{e.fcsF}
\frac{1+s_X}{n_X}+\frac{1+s_Y}{n_Y}+\frac{1+s_Z}{n_Z}\leq F.
\eeq
and
$$\frac{gcd(X,Y)}{rad(gcd(X,Y))}\leq Q,$$ 
\econ
where $rad(g)$ denotes the  product of the distinct prime divisors of $g$.

The finiteness of solutions in the above conjecture is over all possible degrees $n_X,n_Y,n_Z$ 
and spreads satisfying the assumptions.

Note that the above conjecture does not holds when \eqref{e.fcsF} is replaced by \eqref{e.suminverse}.
For example, 
$$(a^\alpha-1)(a^\alpha+1)+1^m=a^{2\alpha}$$
defines an infinite family of examples with $s_X=2,$ $n_X=2,$ and with $n_Y$ and $n_Z$ 
arbitrarily high.
Also, the inequalities \eqref{e.sBound} are necessary, as every number $X$ is a product, $1^n\cdot X$, of an arbitrarily high degree and of spread $X$.

Note that the above conjecture generalizes the Fermat-Catalan conjecture. Indeed, \eqref{e.sBound} holds in this case and
$F$ can be taken to be $41/42$, since it is not difficult to show that $41/42$ is the maximum value of sums $\frac1{n_X}+\frac1{n_Y}+\frac1{n_Z}$ which are less than $1$. More generally, there exists a maximum value for expressions 
\beq\label{e.fcs1}
\frac{s_X+1}{n_X}+\frac{s_X+1}{n_Y}+\frac{s_X+1}{n_Z}<1
\eeq
when $s_X, s_Y$ and $s_Z$ are bounded from the above.

The above conjecture is motivated by Theorem \ref{t.FCabc}, making it a consequence of the abc conjecture.

The Pillai's conjecture states that for any $B>0$ there are finitely many pairs of positive integer powers, which differ by $B$, \cite{ST, Wa, BL}. By restricting our conjecture to $Y=1$ we obtain an analogous statement for products: 

For any $M$ and $F<1$, there are finitely many products $X,Z,$ such that
$$Z-X=B,\quad \frac{s_X^2}{b_X}, \frac{s_Z^2}{b_Z} \leq M,\ \text{and}\ \frac{1+s_X}{n_X}+\frac{1+s_Z}{n_Z}\leq F.$$

\section{A version of the Tijdeman-Zagier conjecture for products}
\label{s.BTZp}

Let us now consider triples $X,Y, Z=X+Y$ of which two terms are powers and, the third one, of the lowest degree is a product of an arbitrary spread. Equivalently, by permuting the variables if necessary, we consider the equation 
$X\pm Y=Z$, where $X,Y$ are powers of degree $n$ and $m$ respectively and $Z$ is a product of degree $d\leq min(m,n).$

On the basis of our computational testing, we propose the following version of Tijdeman-Zagier conjecture for products:

\begin{conjecture}\label{c.GBTZ}
For $x,y$ coprime, $x^n\pm y^m$  is not a product of degree $2<d\leq min(n,m)$ and of spread $s$ satisfying 
\beq\label{e.fcs}
1/n+1/m+ (1+s)/d<1.
\eeq
\end{conjecture}

We discuss the scope of the computational tests of this conjecture in Remark \ref{r.comptest}. 
We do not expect this conjecture 
to be a consequence of the abc conjecture.

\section{
Conjectural Version of Fermat's Last Theorem for products}
\label{s.Fermat}

We call a triple of products $X, Y, Z=X\pm Y$ \underline{F-C} if it satisfies
 \eqref{e.fcs1} and we call it \underline{maxgcd} if $gcd(X,Y)=min(X,Y)$.
It is easy to find maxgcd F-C triples,
for example, 
$$X=x^p,\quad Y=x^{p-1}, \ \text{and}\ Z=x^{p-1}(x+1)$$
for any $p>4$ and any $x.$   
For this reason we will only consider non-maxgcd F-C triples from now on only.
Even among triples of powers there are bountiful non-maxgcd F-C triples. For example, we have:

\bpro\label{p.nonmaxsols}
If $n,m$ are not divisible by $3$ and $gcd(n,m)\leq 2$ then Eq. \eqref{e.FC} has infinitely many non-maxgcd solutions.
\epro

\bpr
For any $a,$ the Pythagorean triple 
$$(a^2-1)^2+(2a)^2=(a^2+1)^2,$$
leads to 
$$(a^2-1)^{2+3m\alpha}(2a)^{3n\beta}(a^2+1)^{nm\gamma}+
(a^2-1)^{3m\alpha}(2a)^{2+3n\beta}(a^2+1)^{nm\gamma}=$$
$$(a^2-1)^{3m\alpha}(2a)^{3n\beta}(a^2+1)^{2+nm\gamma}.$$
By the assumptions on $n$ and $m$, there exist positive $\alpha,\beta,\gamma$  such that
$$n\mid 2+3m\alpha,\quad m\mid 2+3n\beta\quad \ \text{and}\ 3\mid 2+nm\gamma.$$
yielding the desired F-C solutions of \eqref{e.FC}. In fact, there are infinitely many of such pairs $(\alpha,\beta).$
\epr

However, the situation is much different if, as in the previous section, we assume that 
$x^n\pm y^m=Z$ is a product of degree $d\leq min(m,n)$ satisfying the F-C inequality. 
It is natural to expect that the spread of $Z$ is positive for the majority of such F-C triples,  but remarkably the opposite is the case, at least for low values of $d$. 
For example,  there are only four non-maxgcd  F-C solutions 
of $x^n\pm y^m=Z$ of degree $d=3$ and spread $>0$ for $x^n, y^m, Z\leq 2^{80}$:
$$12^4 - 2^{14} = 16^2\cdot 17,\quad 24^4 - 2^{16} = 64^2\cdot 65,$$ 
$$6^8 - 2^{18}=112^2\cdot 113,\ \text{and}\ 117^4 - 3^{14} =567^2\cdot 568.$$
(The largest of the powers and products in these solutions is $117^4<2^{27}.$)

We propose

\bcon\label{c.nonmaxgcd3}
Equation $x^n\pm y^m=Z$ has finitely many non-maxgcd F-C 
 solutions only for $Z$ of degree $d=3$ and spread $s\ne 0$.
\econ

The above may very well be true for higher degrees $d$ as well, however we do not have sufficient computational data to conjecture it.

Our Conjecture \ref{c.GBTZ} can be considered a vast generalization of the statement of the Fermat's Last Theorem.
Based on our computational testing, we propose here another generalization of that theorem, which does not require coprimness of $x$ and $y$:

\bcon\label{c.Fp}
For any $x,y$, the equation $x^n \pm y^n=Z$ has no non-maxgcd F-C solutions
 for products $Z$ of degree $n$.
\econ

(Being an F-C solution means here that the spread of $Z$ is less than $n-3$.)

\brem\label{r.comptest}
We have verified Conjectures \ref{c.GBTZ} and \ref{c.Fp} for values $x^n,y^m,Z \leq M_d,$
where 
$$M_3=2^{80},\ M_4=2^{100}, \ \text{and}\ M_d=2^{113}, \ \text{for}\ 5\leq d\leq 21.$$ 
$d=n$ for Conjecture \ref{c.Fp}. Note that Conjectures \ref{c.GBTZ} and \ref{c.Fp} hold trivially for $d=2$.
We have also tested them for also for higher $d$, but with limited range of $s$ due to the memory requirements and computational complexity of handling products of large spreads. 
\erem
 
As in the case of the conjecture of the previous section, we do not expect the above conjecture to follow from the abc conjecture, 
(although some special cases of it may be consequences of some explicit versions of the abc conjecture).  

More generally, in our testing, we did not find a no non-maxgcd F-C solution to $x^n + y^n=Z$ with a product $Z$ of degree less than $n$.
It is worth mentioning here that \cite{DM, Po} showed that $x^n+y^n=z^2$ has no coprime solutions when $n\geq 4$ and
and  $x^n+y^n=z^3$ has no coprime solutions when $n\geq 3$ assuming modularity of all elliptic curves.
\vspace*{.1in}

\brem
The above conjecture implies that for any $x,y$, a non-max-gcd $x^n \pm y^n=Z$ cannot be a product of $d\leq n$ numbers, each less than $n-2.$ Indeed, if it is then multiplied by $1^{n-d}$ then it results in a product of degree $n$ and spread $s<n-3.$
\erem

More generally, the possible combinations $(n,m,d)$ for non-maxgcd F-C solutions of the equation $x^n\pm y^m=Z$ 
are quite restricted, going beyond the above conjectures. For example, as mentioned above, we did not find any such solution of
$x^n\pm y^n=Z$ with $Z$ of degree $d<n.$ 

\bque
For what values $d,n,m$, does the equation $x^n\pm y^m=X$ have a non-maxgcd F-C solution?
\eque

Let us comment more about maxgcd solutions now.
For any $n$ there is an infinite family of them
 $$(v\cdot w^n)^n+(v\cdot w^{n-1})^n=(v\cdot w^n)^{n-1}(w^n+v)$$
(for any $v,w$). We call such solutions \underline{standard}.

\bcon\label{c.main} 
All maxgcd F-C solutions of $x^n\pm y^n=Z$ for $Z$ of degree $n$ and spread $1$ are standard.
\econ

\section{Relation to the abc conjecture.}
\label{s.abc}

The celebrated abc conjecture of Oesterl\'e–Masser, \cite{Oe}, states
\vspace*{.1in}

\noindent{\bf abc Conjecture}
{\it For every $\ve>0$ there is $C_\ve$ such that for all positive coprime pairs of integers $a,b$
$$a+b< C_\ve\cdot rad(ab(a+b))^{1+\ve}.$$}
\vspace*{.1in}

Shinichi Mochizuki developed inter-universal Teichmüller theory and asserted that it implies the above statement, see in particular \cite{Mo}. The correctness of his theory and completeness of the proof of the abc conjecture have not been however widely accepted in Number Theory community; see more detailed discussion in \cite{Wik}.

The abc conjecture has far reaching consequences for number theory, see eg. \cite{El, GT}. In particular,  it implies the Fermat-Catalan conjecture. However, its power is tempered by the non-explicity of the constants $C_\ve.$ Certain explicit versions of the abc conjecture were proposed, eg. in \cite{Ba,Go,Wa}. 
Although not directly related to the main topic of this paper, our explicit version of the abc conjecture may be of independent interest to the reader:

\bcon\label{c.abc-mine}
For all positive coprime pairs of integers $a,b$
$$a+b =c< max(rad(ab), rad(ac), rad(bc))\cdot rad(abc)^{7/8}.$$
\econ

We confirmed the correctness of this conjecture for all triples consisting of $a,b,c< 2^{63}$, based on the data available at \cite{ABC}. 

Although the standard formulation of the abc conjecture concerns coprime pairs $(a,b),$ it will be convenient for us to drop this assumption. 

\bpro
\label{p.abc2}
Assuming the abc conjecture, for every $Q$ and $\ve>0$ there are only finitely many positive integers $a,b$ such that 
\beq\label{e.abcin}
a+b>rad(ab(a+b))^{1+\ve}
\eeq
and
\beq\label{e.gcdin}
\frac{gcd(a,b)}{rad(gcd(a,b))}\leq Q.
\eeq
\epro

\bpr
Assume the abc conjecture and assume, to the contrary, that there are infinitely many pairs $(a,b)$ satisfying the above inequalities.
For any such pair, let $g=gcd(a,b)$ and $a'=a/g, b'=b/g$ and $c'=a'+b'.$ Then we have 
$$c'g>rad(a'b'c'g)^{1+\ve}=R^{1+\ve}rad(g)^{1+\ve},$$
for them, where $R=rad(a'b'c')$.
But since 
\beq\label{e.c'}
c'<R^{1+\ve/2} 
\eeq 
for almost all of them, we have 
$$R^{1+\ve/2} g> R^{1+\ve}rad(g)^{1+\ve}$$ 
for an infinity of them, i.e.
$$R^{1+\ve}/R^{1+\ve/2}< g/rad(g)^{1+\ve}<Q.$$
Hence, $R$ is bounded from above, $c'$ is bounded from above as well, by \eqref{e.c'}. 
However, that is possible for finitely many $a',b'$ only.

This implies that among the infinity of pairs $(a,b)$ as above, there is an infinite subfamily 
with $a/g=a'$ and $b/g=b'$ for certain fixed $a'$ and $b',$ where $g=gcd(a,b).$
For these pairs, inequality \eqref{e.abcin}
states
$$(a'+b')g>rad(a'b'(a'+b'))^{1+\ve}g^{1+\ve}.$$
Therefore, the values of $g$ must be bounded from above, yielding a contradiction.
\epr

The statement of Proposition \ref{p.abc2}, in turn, even for $Q=1$, implies the abc conjecture, cf. \cite{MM}.

\bthm\label{t.FCabc}
The abc conjecture implies Conjecture \ref{c.FCS}.
\ethm

As before, given a product $X$ we denote its degree by $d_X$, its span by $s_X$, and we assume that $X$ is a product of $d_X$ numbers from $b_X,...,b_X+s_X$ (possibly with multiplicities). 
We require that $s_X+1<d_X.$
Let us start with the following:

\blem
For any product $X$, 
we have 
$$rad\, X\leq e^{2s_X^2/b_X}\cdot X^{(s_X+1)/d_X}.$$
\elem

\bpr
Since $1+s_X/b_X\leq e^{s_X/b_X}$, we have $b_X+s_X\leq e^{s_X/b_X}b_X$ and
$$(b_X+s_X)^{d_X}\leq e^{s_X(d_X-1)/b_X}b_X^{d_X-1}(b_X+s_X)\leq e^{s_X(d_X-1)/b_X} X.$$
Hence,
$$rad\, X\leq (b_X+s_X)^{s_X+1}\leq e^{s_X(s_X+1)(d_X-1)/(b_Xd_X)}X^{(s_X+1)/d_X}.$$
Since 
$$(s_X+1)(d_X-1)/d_X\leq 2s_X,$$ 
the statement follows.
\epr

\noindent{\it Proof of Theorem \ref{t.FCabc}:} Let us prove finiteness of the number of Fermat-Catalan solutions of 
$A\cdot X+B\cdot Y=C\cdot Z$ 
satisfying \eqref{e.gcdin}.
By Proposition \ref{p.abc2}, there are only finitely many such triples satisfying
$$CZ> rad(ABXYZ)^{1+\ve},$$
for any $\ve>0.$
Since each such $Z$ corresponds to finitely many solutions only, it is enough to show that there is only finitely many solutions $(X,Y,Z)$
satisfying 
\beq\label{e.abc}
CZ\leq  rad(ABXYZ)^{1+\ve}
\eeq
for some $\ve<1/F-1$.

By the above lemma, we have
$$rad(XYZ)\leq e^{6M}\cdot X^{(s_X+1)/d_X}Y^{(s_Y+1)/d_Y}Z^{(s_Z+1)/d_Z}\leq e^{6M}\cdot Z^f,$$
where $f=(s_X+1)/d_X + (s_Y+1)/d_Y +(s_Z+1)/d_Z\leq F<1.$

Hence, \eqref{e.abc} implies 
$$Z\leq \frac{AB}C(e^{6M} Z^f)^{1+\ve}.$$
Since $f(1+\ve)\leq F<1$ 
inequality  \eqref{e.abc} has finitely many solutions $Z$ only.
\qed

%

\end{document}